\documentclass{amsart}

\usepackage{amssymb}
\usepackage{amsthm}
\usepackage{graphicx,epsfig}
\usepackage{etoolbox}
\usepackage{url}

\newcommand{\R}{\mathbb{R}}

\newcommand{\C}{\mathbb{C}}
\newcommand{\Z}{\mathbb{Z}}

\newcommand{\Znz}{\Z_{\neq 0}}
\newcommand{\N}{\Z_{>0}}
\newcommand{\bigoh}{{\mathcal{O}}}
\newcommand{\e}{{\rm e}}

\newcommand{\D}{{\rm d}}
\newcommand{\erfc}{{\rm erfc}}
\newcommand{\Phihat}{{\widehat{\Phi}}}
\newcommand{\phihat}{{\widehat{\phi}}}
\newcommand{\Ei}{{\rm Ei}}

\newcommand{\erf}{{\rm erf}}

\newtoggle{longver}
\togglefalse{longver}

\newtheorem{theorem}{Theorem}[section]
\newtheorem{lemma}[theorem]{Lemma}

\numberwithin{equation}{section}

\begin{document}

\title{Computing $\pi\left(x\right)$ Analytically}


\author{David J. Platt}
\address{Heilbronn Institute for Mathematical Research, University of Bristol}
\email{dave.platt@bris.ac.uk}
\thanks{\\
\indent{Accepted for publication in Math. Comp. 8 Sept. 2013}\\
\\
{\hangindent=\parindent This work formed part of my PhD research and I would like to thank Dr.\ Andrew Booker for his patient supervision. Funding was provided by the Engineering and Physical Sciences Research Council through the University of Bristol and I am grateful to both.}}



\date{}


\begin{abstract}
We describe a rigorous implementation of the Lagarias and Odlyzko Analytic
Method to evaluate the prime counting function and its use to compute
unconditionally the number of primes less than $10^{24}$.
\iftoggle{longver}{This includes a new, FFT based algorithm for computing Riemann's $\zeta$ function on the half line.}{}
\end{abstract}

\maketitle

\section{Introduction}

Computing exact values of the function $\pi(x)$, which counts the number of primes less than
or equal to $x$, has exercised mathematicians since antiquity. Early methods
involved enumerating all the primes less than the target $x$ (using, for
example, the sieve of Eratosthenes) and then counting them.  In $1870$ Meissel \cite{Meissel1870} described a combinatorial method which he eventually used to manually compute $\pi\left(10^9\right)$ \cite{Meissel1885} (albeit not quite accurately). The algorithm was subsequently improved by Lehmer \cite{Lehmer1959}, then by Lagarias, Miller and Odlyzko \cite{Lagarias1985} and most recently by Del\'{e}glise and Rivat \cite{Deleglise1996}. In $2007$ Oliveira e Silva used the algorithm to compute $\pi\left(10^{23}\right)$.

The Prime Number Theorem dictates that all methods reliant on enumerating the primes must have time complexity of $\Omega\left(x\log^{-1} x\right)$. The latest incarnations of the combinatorial method achieve $\bigoh\left(x^{2/3}\log^{-2} x\right)$.

In their 1987 paper \cite{Lagarias1987}, Lagarias and Odlyzko described an
analytic algorithm with (in one form) time complexity
$\bigoh\left(x^{1/2+\epsilon}\right)$. In $2010$ B\"uthe, Franke, Jost and
Kleinjung announced a value for $\pi\left(10^{24}\right)$ \cite{Buthe2013a}
contingent on the Riemann Hypothesis. Their approach ``is similar to the one
described by Lagarias and Odlyzko, but uses the Weil explicit formula instead
of complex curve integrals''. This paper describes an implementation reverting to Riemann's explicit formula which we have used to compute
$\pi\left(10^{24}\right)$ unconditionally.

\section{A Note on Rigour}

To many, rigorous computation is an oxymoron, due to potential bugs
in hardware, operating systems, compilers and (most likely) user's code. Add to this the
chance that a power spike or cosmic ray interaction could scupper even a
correctly written application, and the situation seems hopeless. We do what we
reasonably can to minimise such risks including running applications on
systems with ECC memory after testing on hardware from different vendors with
different
  operating systems and using different compilers.

However, there are certain aspects over which we do have more
control. Estimating the rounding error that will accumulate through a complex
floating point computation is a non-trivial task that we eschew. Rather, we rely on interval
arithmetic (see \cite{Moore1966} for a good introduction). Thus, instead of storing a single floating point
approximation, we hold an interval comprising two floating points to bracket the true
value. We then overload the standard operators and functions to
handle such intervals. 

Furthermore, we will estimate some quantities by, for example,
truncating an infinite sum. We will need to derive a rigorous bound for the
error introduced, but rather than keep track of such errors manually, we can
simply add them to the interval being evaluated and let the software take
over. Indeed, we can in some circumstances use interval arithmetic to
compute such bounds for us.

\section{The Analytic Algorithm}

The analytic algorithm relies on Perron's formula.

\begin{theorem}[Perron's Formula]
Let $a(n)$ be an arithmetic function with Dirichlet series
\begin{equation*}
g(s)=\sum\limits_{n=1}^\infty \frac{a(n)}{n^s}.
\end{equation*}
Now if $g(s)$ is absolutely convergent whenever $\Re s>\sigma_a$, then for $c>\sigma_a$ and $x>0$ we have
\begin{equation*}
\sum\limits_{n\leq x}{^*}a(n)=\frac{1}{2\pi i}\int\limits_{c-i\infty}^{c+i\infty}g(s)x^s \frac{\D s}{s}
\end{equation*}
where the $^*$ on the summation sign indicates that if $x$ is an integer then
only $1/2$ of the $a(x)$ term is included.
\begin{proof}
See page 245 of \cite{Apostol1976} and the subsequent note.
\end{proof}
\end{theorem}

The relevance of Perron's formula to the matter at hand comes from the series,
absolutely convergent for $\Re s>1$, 
\begin{equation*}
\log\zeta(s)=\sum\limits_{n=2}^\infty\frac{\Lambda(n)}{n^s\log n}.
\end{equation*}
Here $\Lambda$ is the von-Mangoldt function so $\frac{\Lambda(n)}{\log n}$
is $\frac{1}{m}$ at prime powers $p^m$ and zero elsewhere. Define
\begin{equation*}
\pi^*(x):=\sum\limits_{p^m\leq x}\frac{1}{m}
\end{equation*}
where if $x$ is a prime power we only take $1/2$ of its contribution
to the sum. Then applying Perron's formula we get
for $c>1$ and $x>0$
\begin{equation}\label{eq:pi_star}
\sum\limits_{n\leq x}{^*}\frac{\Lambda(n)}{\log n}=\pi^*(x)=\frac{1}{2\pi i}\int\limits_{c-i\infty}^{c+i\infty}\log\zeta(s)x^s \frac{\D s}{s}.
\end{equation}
Here we note that, although we can cheaply recover $\pi(x)$ from $\pi^*(x)$,
the slow rate of convergence of the integral dooms any attempt to use it in
this context.

At this point, Lagarias and Odlyzko introduce a ``suitable'' Mellin transform pair $\phi(t)$ and $\hat\phi(s)$ and derive
\begin{equation}\label{eq:anal_pi}
\pi^*(x)=\frac{1}{2\pi i}\int\limits_{\sigma-i\infty}^{\sigma+i\infty}\log\zeta(s)\hat\phi(s)\D s+\sum\limits_{p^m}\frac{1}{m}\left[\chi_x(p^m)-\phi(p^m)\right]
\end{equation}
where $\chi_x(t)$ is defined by
\begin{equation}
\nonumber
\chi_x(t):=\left\{\begin{array}{ll}1\;\;\;&t<x\\1/2&t=x.\\0&t>x\end{array}\right.
\end{equation}
We note that taking $\hat\phi(s)=\frac{x^s}{s}$ makes $\phi(t)=\chi_x(t)$ and
we recover (\ref{eq:pi_star}).

Thus estimating $\pi^*(x)$ now splits into estimating an integral and summing
$\phi(t)$ evaluated at prime powers in the vicinity of $x$.

In his PhD thesis \cite{Galway2004} Galway investigated the proposed algorithm and suggested using the Mellin transform pair
\begin{equation}\label{eq:phihat}
\hat\phi(s):=\frac{x^s}{s}\exp\left(\frac{\lambda^2s^2}{2}\right)\textrm{ and } \phi(t):=\frac{1}{2}\erfc\left(\frac{\log\left(\frac{t}{x}\right)}{\sqrt{2}\lambda}\right).
\end{equation}
Here $\erfc$ is the complementary error function
\begin{equation}
\nonumber
\erfc(x):=\frac{2}{\sqrt\pi}\int\limits_x^\infty\exp\left(-t^2\right)\D t
\end{equation}
and $\lambda$ is a positive real parameter used to balance the convergence of
the integral with the width of the prime sieve.

Galway showed that $\phi$ and $\hat\phi$ as defined in (\ref{eq:phihat}) are indeed ``suitable''
and, using arguments based on the uncertainty principle, suggested that they
are in some sense optimal. He also gave a rigorous bound for the
error introduced by truncating the prime sieve to some finite width.



\section{Evaluating $\frac{1}{2\pi
    i}\int\limits_{\sigma-i\infty}^{\sigma+i\infty}\log\zeta(s)\phihat(s)\D s$}

At this point, we depart from the line taken by Lagarias and Odlyzko. Rather
than attempt to numerically estimate the integral in (\ref{eq:anal_pi}), we take an approach closer to the spirit of Riemann and evaluate it in terms of the non-trivial zeros of $\zeta$, leading to Theorem \ref{th:Gsum} below.

Before proceeding, we need a couple of lemmas.
\begin{lemma}\label{lem:rtp} The ``Round the Pole'' lemma. Let $f$ be a meromorphic function with a simple pole at $\alpha$ with residue R, and let $\Gamma$ be the semicircular contour anticlockwise from $\alpha+\epsilon$ to $\alpha-\epsilon$. Then
\begin{equation}
\nonumber
\lim\limits_{\epsilon\rightarrow 0^+}\int\limits_\Gamma f(z)\D z=\pi i R.
\end{equation}
\begin{proof}
See page 29 of \cite{Open1994}.
\end{proof}
\end{lemma}

\begin{lemma}\label{lem:logzetaeps}
We have
\begin{equation}
\nonumber
\lim\limits_{\epsilon\rightarrow 0^+}(\log\zeta(1+\epsilon)-\log(-\zeta(1-\epsilon)))=0.
\end{equation}
\begin{proof}
\begin{equation}
\nonumber
\begin{aligned}[l]
\lim\limits_{\epsilon\rightarrow 0^+}(\log\zeta(1+\epsilon)-\log(-\zeta(1-\epsilon)))&=\lim\limits_{\epsilon\rightarrow 0^+}\log\frac{\zeta(1+\epsilon)}{-\zeta(1-\epsilon)}\\
&=\lim\limits_{\epsilon\rightarrow 0^+}\log\frac{1/\epsilon+\bigoh(1)}{1/\epsilon+\bigoh(1)}=0.\hfill\qedhere
\end{aligned}
\end{equation}
\end{proof}
\end{lemma}

\begin{lemma}\label{lem:zetapbound}
There exists a sequence of $T_j\rightarrow\infty$ such that for any $\sigma\in[-1,2]$ we have for $s=\sigma+iT_j$
\begin{equation}
\nonumber
\frac{\zeta '}{\zeta}(s)=\bigoh(\log^2 T_j).
\end{equation}
\begin{proof}

Referring to Davenport \cite{Davenport2000}, for any zero $\beta+i\gamma$ of $\zeta$ with $\gamma$ large, we note that there are $\bigoh(\log \gamma)$ zeros with imaginary part $\in[\gamma-1,\gamma+1]$ (Corollary (a), page $99$). Therefore we can select a $T_j$ within $\bigoh(1)$ of $\gamma$ such that $T_j$ differs from the imaginary part of any zero by $\gg \frac{1}{\log T_j}$. By (4) on page $99$ we have for $\sigma\in[-1,2]$
\begin{equation}
\nonumber
\left|\frac{\zeta '}{\zeta}\left(\sigma+iT_j\right)\right|=\left|\sum\limits_\rho\, ^{'} \frac{1}{\sigma+iT_j-\rho} +\bigoh(\log T_j)\right|
\end{equation}
where the sum is taken over zeros with imaginary part $\in[T_j-1,T_j+1]$. There are $\bigoh(\log T_j)$ such zeros, each one making a contribution to the sum limited by $\bigoh(\log T_j)$ and the result follows.
\end{proof}
\end{lemma}

\begin{lemma}\label{lem:log_zeta_bound}
For $t\in\R$
\begin{equation}
\nonumber
\left|\log(-\zeta(-1+it))\right|\leq 5+t^2.
\end{equation}
\begin{proof}
By the functional equation
\begin{equation}
\nonumber
\zeta(-1+it)=\zeta(2-it)\frac{\Gamma\left(\frac{2-it}{2}\right)}{\Gamma\left(\frac{-1+it}{2}\right)}\pi^{\left(-\frac{3}{2}+it\right)}.
\end{equation}
We then use $\left(\frac{-1+it}{2}\right)\Gamma\left(\frac{-1+it}{2}\right)=\Gamma\left(\frac{1+it}{2}\right)$ so we can apply Stirling's approximation. Also, for $\Re s>1$ we have
\begin{equation*}
\begin{aligned}
\left|\log\zeta(s)\right|=\left|\sum\limits_{n=2}^\infty\frac{\Lambda(n)}{\log(n)} \frac{1}{n^s}\right|\leq \sum\limits_{n=2}^\infty\frac{\Lambda(n)}{\log(n)} \frac{1}{n^{\Re s}}= \log\zeta(\Re s).&\qedhere
\end{aligned}
\end{equation*}
\end{proof}
\end{lemma}

\begin{lemma}\label{lem:minus_one_int}
\begin{equation}
\nonumber
\left|\frac{1}{2\pi i}\int\limits_{-1-i\infty}^{-1+i\infty}\log(-\zeta(s))\phihat(s)\D s\right| \leq\frac{\exp\left(\frac{\lambda^2}{2}\right)}{2\pi x\lambda}\left(5\sqrt{2\pi}+\frac{2}{\lambda}\right).
\end{equation}
\begin{proof}
We use Lemma \ref{lem:log_zeta_bound} and take absolute values, majorising with
\begin{equation}
\nonumber
\frac{\exp\left(\frac{\lambda^2}{2}\right)}{2\pi x}\left[5\int\limits_{-\infty}^\infty\exp\left(\frac{-\lambda^2 t^2}{2}\right)\D t +\int\limits_{-\infty}^\infty \left|t\right|\exp\left(\frac{-\lambda^2 t^2}{2}\right)\D t\right]
\end{equation}
where both integrals can be evaluated.
\end{proof}
\end{lemma}

\begin{lemma}\label{lem:Phihatdef}
Let $\Phihat(s)$ be the unique holomorphic function $\Phihat:\C\setminus\R_{\leq 0}\rightarrow\C$ such that
\begin{itemize}
\item $\Phihat'(s)=\phihat(s)$ and
\item $\lim\limits_{t\rightarrow\infty}\left[\Phihat(\sigma+it)+\Phihat(\sigma-it)\right]=0$ for any fixed real $\sigma$.
\end{itemize}
Then
\begin{enumerate}
\item $\Phihat(s)-\log s$ extends to an entire function,
\item $\lim\limits_{t\rightarrow\infty}\Phihat(\sigma+it)=C$ is purely
  imaginary and
\item $\Phihat(\sigma\pm it)\mp C$ is rapidly decreasing as $t\rightarrow\infty$.
\end{enumerate}
\begin{proof}
To show (1) we define for $s\not\in\R_\leq{0}$
\begin{equation}
\nonumber
F(s)=\int\limits_1^s\phihat(z)\D z
\end{equation}
where the contour of integration is the straight line from $1$ to $s$. We then define
\begin{equation}
\nonumber
\Phihat(s):=\lim\limits_{T\rightarrow\infty}\left[F(s)-\frac{F(1+iT)+F(1-iT)}{2}\right]
\end{equation}
and we have
\begin{equation}
\nonumber
F(s)-\log s=\int\limits_1^s\left(\phihat(z)-\frac{1}{z}\right)\D z.
\end{equation}
To show (2) we observe that $\Phihat(\overline{s})=\overline{\Phihat(s)}$ and from the definition we have $C+\overline{C}=0$.
To show (3), we take $T>0$ and we have
\begin{equation*}
\begin{aligned}
|\Phihat(\sigma\pm iT)\mp C|\leq \frac{x^\sigma}{T}\exp\left(\frac{\sigma^2\lambda^2}{2}\right)\int\limits_T^\infty \exp\left(\frac{-\lambda^2t^2}{2}\right)\D t.&\qedhere
\end{aligned}
\end{equation*}
\end{proof}
\end{lemma}



\begin{theorem}\label{th:Gsum}
Let $\Phihat(s)$ be defined as in Lemma \ref{lem:Phihatdef}. Then
\begin{equation}
\nonumber 
\frac{1}{2\pi i}\int\limits_{2-i\infty}^{2+i\infty}\phihat(s)\log\zeta(s)\D
s=\Phihat(1)-\sum\limits_{\rho}\Re \Phihat(\rho)-\log(2)+\frac{1}{2\pi i}\int\limits_{-1-i\infty}^{-1+i\infty}\phihat(s)\log(-\zeta(s))\D s.
\end{equation}
\begin{proof}

We will refer to the contours represented in Figure \ref{fig:contour}.
\begin{figure}[ptb]
\centering
\epsfig{file=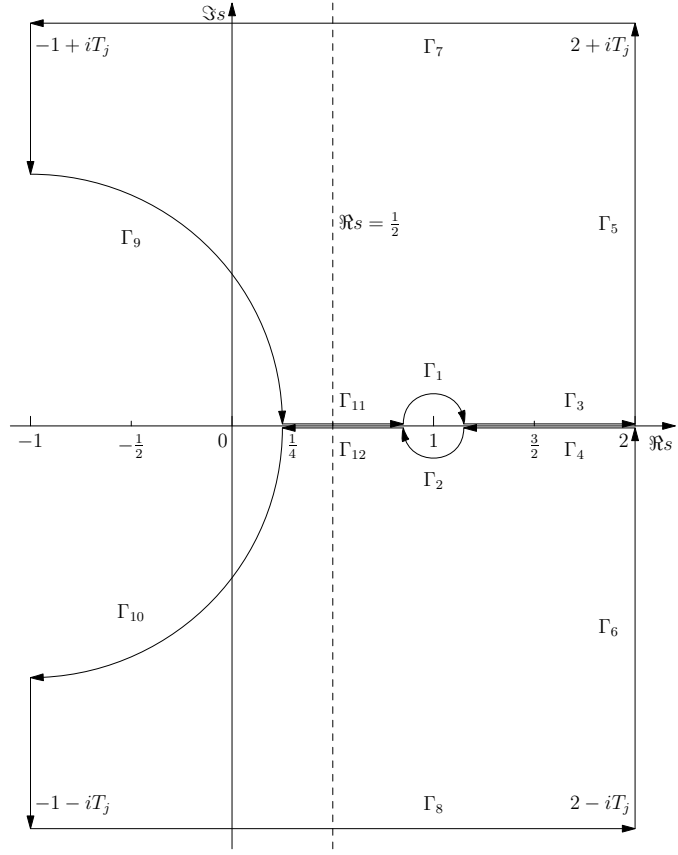,width=0.7\linewidth,clip=}
\caption{Contours to evaluate $\frac{1}{2\pi i}\int\limits_{2-i\infty}^{2+i\infty} \hat\phi(s)\log\zeta(s)\D s$}
\label{fig:contour}
\end{figure}
These contours are
\begin{itemize}
\item $\Gamma_1$ - the semi-circle clockwise from $1-\epsilon$ to $1+\epsilon$ for $\epsilon$ small and positive.
\item $\Gamma_2$ - the semi-circle clockwise from $1+\epsilon$ to $1-\epsilon$.
\item $\Gamma_3$ - the horizontal line from $1+\epsilon$ to $2$.
\item $\Gamma_4$ - the horizontal line from $2$ to $1+\epsilon$.
\item $\Gamma_5$ - the vertical line from $2$ to $2+iT_j$, $T_j$ not the ordinate of a zero of $\zeta$.
\item $\Gamma_6$ - the vertical line from $2-iT_j$ to $2$.
\item $\Gamma_7$ - the horizontal line from $2+iT_j$ to $-1+iT_j$.
\item $\Gamma_8$ - the horizontal line from $-1-iT_j$ to $2-iT_j$.
\item $\Gamma_9$ - the vertical line from $-1+iT_j$ to $-1+\frac{5}{4}i$, followed by the clockwise circular arc centred at $-1$ to $\frac{1}{4}$.
\item $\Gamma_{10}$ - the clockwise circular arc centred at $-1$ from $\frac{1}{4}$ to $-1-\frac{5}{4}i$, followed by the vertical line to $-1-iT_j$. 
\item $\Gamma_{11}$ - the horizontal line from $\frac{1}{4}$ to $1-\epsilon$.
\item $\Gamma_{12}$ - the horizontal line from $1-\epsilon$ to $\frac{1}{4}$.
\end{itemize}

We consider the integrals
\begin{equation}\label{eq:int_top}
\frac{1}{2\pi i}\int(\Phihat(s)-C)\frac{\zeta'(s)}{\zeta(s)}\D s
\end{equation}
for the contours $\Gamma_1$, $\Gamma_3$, $\Gamma_5$, $\Gamma_7$, $\Gamma_9$ and $\Gamma_{11}$ in the upper half plane and
\begin{equation}\label{eq:int_bot}
\frac{1}{2\pi i}\int(\Phihat(s)+C)\frac{\zeta'(s)}{\zeta(s)}\D s
\end{equation}
for $\Gamma_2$, $\Gamma_4$, $\Gamma_6$, $\Gamma_8$, $\Gamma_{10}$ and $\Gamma_{12}$ in the lower half. 

We denote the integrals in (\ref{eq:int_top}) or (\ref{eq:int_bot}) as appropriate along $\Gamma_n$ by $I_n$ and proceed as follows.

For $I_5$ and $I_6$ we get
\begin{equation}
\nonumber
\begin{aligned}
\lim\limits_{j\rightarrow\infty}(I_5+I_6)&=\lim\limits_{j\rightarrow\infty}\frac{1}{2\pi i}\left[\int\limits_{\Gamma_5}\left(\Phihat(s)-C\right)\frac{\zeta '}{\zeta}(s)\D s+\int\limits_{\Gamma_6}\left(\Phihat(s)+C\right)\frac{\zeta '}{\zeta}(s)\D s\right]\\
&=\lim\limits_{j\rightarrow\infty}\frac{1}{2\pi i}\left[\left.\left(\Phihat(s)-C\right)\log\zeta(s)\right|_2^{2+iT_j}+\left.\left(\Phihat(s)+C\right)\log\zeta(s)\right|_{2-iT_j}^2\right]\\
&\;\;\;\;\;\;\;-\frac{1}{2\pi i}\left[\int\limits_{\Gamma_{5,6}}\phihat(s)\log\zeta(s)\D s\right]\\
&=\frac{1}{2\pi i}\left[ 2C\log\zeta(2)-\int\limits_{2-i\infty}^{2+i\infty}\phihat(s)\log\zeta(s)\D s\right]
\end{aligned}
\end{equation}
where $\Gamma_{5,6}$ denotes the contour $\Gamma_5$ followed by $\Gamma_6$. 

Considering the contours $\Gamma_7$ and $\Gamma_8$, we use Lemma \ref{lem:zetapbound} and the Gaussian decay of $\Phihat(s)\pm C$ from Lemma \ref{lem:Phihatdef} to conclude

\begin{equation}
\nonumber
\begin{aligned}
\lim\limits_{j\rightarrow\infty}(I_7+I_8)&=\lim\limits_{j\rightarrow\infty}\frac{1}{2\pi i}\left[\int\limits_{\Gamma_7}\left(\Phihat(s)-C\right)\frac{\zeta '}{\zeta}(s)\D s+\int\limits_{\Gamma_8}\left(\Phihat(s)+C\right)\frac{\zeta '}{\zeta}(s)\D s\right]\\
&=0.
\end{aligned}
\end{equation}

Considering $I_9$ and $I_{10}$ we have
\begin{equation}
\nonumber
\begin{aligned}
\lim\limits_{j\rightarrow\infty}(I_{9}+I_{10})&=\lim\limits_{j\rightarrow\infty}\frac{1}{2\pi i}\left[\int\limits_{\Gamma_{9}}\left(\Phihat(s)-C\right)\frac{\zeta '}{\zeta}(s)\D s+\int\limits_{\Gamma_{10}}\left(\Phihat(s)+C\right)\frac{\zeta '}{\zeta}(s)\D s\right]\\
&=\lim\limits_{j\rightarrow\infty}\frac{1}{2\pi i}\left[\left.\left(\Phihat(s)-C\right)\log(-\zeta(s))\right|^{1/4}_{-1+iT_j}\right.\\
&\;\;\;\;\;\;\;\;\;\;\;\;\;\;+\left.\left.\left(\Phihat(s)+C\right)\log(-\zeta(s))\right|_{1/4}^{-1-iT_j}\right]\\
&\;\;\;\;\;\;\;-\frac{1}{2\pi i}\left[\int\limits_{\Gamma_9}\phihat(s)\log(-\zeta(s))\D s+\int\limits_{\Gamma_{10}}\phihat(s)\log(-\zeta(s))\D s\right]\\
&=-\frac{1}{2\pi i}\left[\int\limits_{\Gamma_9,\Gamma_{10}}\phihat(s)\log(-\zeta(s))\D s+2C\log\left(-\zeta\left(1/4\right)\right)\right]
\end{aligned}
\end{equation}
where the contour of integration is $\Gamma_9$ followed by $\Gamma_{10}$. Convergence of this integral is due to Lemma
\ref{lem:minus_one_int} and the zero free region of $\zeta(s)$ with
$|s+1|\leq\frac{5}{4}$ and $\Re s \geq -1$.

For $I_{11}$ and $I_{12}$ we have
\begin{equation}
\nonumber
\begin{aligned}
I_{11}+I_{12}&=\frac{1}{2\pi i}\left[\int\limits_{\Gamma_{11}}\left(\Phihat(s)-C\right)\frac{\zeta '}{\zeta}(s)\D s+\int\limits_{\Gamma_{12}}\left(\Phihat(s)+C\right)\frac{\zeta '}{\zeta}(s)\D s\right]\\
&=\frac{1}{2\pi i}\left[\left.\left(\Phihat(s)-C\right)\log(-\zeta(s))\right|^{1-\epsilon}_{1/4}+\left.\left(\Phihat(s)+C\right)\log(-\zeta(s))\right|_{1-\epsilon}^{1/4}\right]\\
&\;\;\;\;\;\;\;-\frac{1}{2\pi i}\left[\int\limits_{\Gamma_{11}}\phihat(s)\log(-\zeta(s))\D s+\int\limits_{\Gamma_{12}}\phihat(s)\log(-\zeta(s))\D s\right]\\
&=\frac{1}{2\pi i}\left[2C\log\left(-\zeta\left(1/4\right)\right)-2C\log(-\zeta(1-\epsilon))\right].
\end{aligned}
\end{equation}

For $I_1$ and $I_2$ we find
\begin{equation}
\nonumber
\begin{aligned}
I_{1}+I_{2}&=\frac{1}{2\pi i}\left[\int\limits_{\Gamma_1}\left(\Phihat(s)-C\right)\frac{\zeta '}{\zeta}(s)\D s+\int\limits_{\Gamma_2}\left(\Phihat(s)+C\right)\frac{\zeta '}{\zeta}(s)\D s\right]\\
&=\frac{1}{2\pi i}\left[\oint\Phihat(s)\frac{\zeta'(s)}{\zeta(s)}\D s-C\int\limits_{\Gamma_1}\frac{\zeta '}{\zeta}(s)\D s+C\int\limits_{\Gamma_2}\frac{\zeta '}{\zeta}(s)\D s\right]\\
&=\Phihat(1)-\frac{C}{2\pi i}\left[\int\limits_{\Gamma_1}\frac{\zeta '}{\zeta}(s)\D s-\int\limits_{\Gamma_2}\frac{\zeta '}{\zeta}(s)\D s\right]
\end{aligned}
\end{equation}
by Cauchy's Theorem since the residue of $\frac{\zeta'}{\zeta}(s)$ at $s=1$ is $-1$.

Finally, for $I_3$ and $I_4$ we get
\begin{equation}
\nonumber
\begin{aligned}
I_3+I_4&=\frac{1}{2\pi i}\left[\int\limits_{\Gamma_3}\left(\Phihat(s)-C\right)\frac{\zeta '}{\zeta}(s)\D s+\int\limits_{\Gamma_4}\left(\Phihat(s)+C\right)\frac{\zeta '}{\zeta}(s)\D s\right]\\
&=\frac{1}{2\pi i}\left[\left.\left(\Phihat(s)-C\right)\log\zeta(s)\right|_{1+\epsilon}^2+\left.\left(\Phihat(s)+C\right)\log\zeta(s)\right|^{1+\epsilon}_2\right]\\
&\;\;\;\;\;\;\;-\frac{1}{2\pi i}\left[\int\limits_{\Gamma_3}\phihat(s)\log\zeta(s)\D s+\int\limits_{\Gamma_4}\phihat(s)\log\zeta(s)\D s\right]\\
&=\frac{1}{2\pi i}\left[2C\log\zeta(1+\epsilon)-2C\log\zeta(2)\right].
\end{aligned}
\end{equation}

Now by Cauchy's Theorem and exploiting the fact that the non-trivial zeros of $\zeta$ occur in complex conjugate pairs, $\lim\limits_{j\rightarrow\infty}\sum_{k=1}^{12}I_k=\sum\limits_{\rho} \Re \Phihat(\rho)$ so we have

\begin{equation}
\nonumber
\begin{aligned}
\sum\limits_{\rho} \Re \Phihat(\rho)&=\Phihat(1)-\frac{1}{2\pi i}\left[\int\limits_{2-i\infty}^{2+i\infty}\phihat(s)\log\zeta(s)\D
s+\int\limits_{\Gamma_9,\Gamma_{10}}\phihat(s)\log(-\zeta(s))\D
s\right]\\
&+\frac{C}{\pi i}\left[\log \zeta(1+\epsilon)-\log(-\zeta(1-\epsilon))\right]\\
&-\frac{C}{2\pi i}\left[\int\limits_{\Gamma_1}\frac{\zeta'(s)}{\zeta(s)}\D s -\int\limits_{\Gamma_2}\frac{\zeta'(s)}{\zeta(s)}\D s\right].
\end{aligned}
\end{equation}
Now the result follows from taking the limit as $\epsilon\rightarrow 0^+$ by Lemmas
\ref{lem:logzetaeps} and \ref{lem:rtp} and then straightening the line of integration of the second integral to $\Re s=-1$. This introduces a contribution of $\log(-\zeta\left(0\right))=-\log 2$ from the pole of $\phihat(s)$ at $s=0$ with residue $1$.
\end{proof}
\end{theorem}

Again, if we take $\phihat(s)=\frac{x^s}{s}$ then $\Phihat(s)=\Ei(\log s)$
where $\Ei$ is the exponential integral and we recover Riemann's explicit
formula
\begin{equation*}
\pi^*(x)=\Ei(\log x)-\sum\limits_\rho \Ei(\rho\log x)-\log 2
+\int\limits_{-1-i\infty}^{-1+i\infty}\log(-\zeta(s))x^s\frac{\D s}{s}.
\end{equation*}

We truncate the sum over zeros so we need a rigorous bound for the error that
this introduces. We derive such a bound in Appendix \ref{app:zero_trunc}.

The computation of $\pi(x)$ now reduces to
\begin{itemize}
\item enumerating the prime powers near $x$,
\item computing $\phi(t)$ at these prime powers,
\item locating the non-trivial zeros of $\zeta$ to sufficient accuracy and
\item evaluating $\Phihat$ at these zeros (and at $1$).
\end{itemize}

\section{The Prime Sieve and $\phi(p)$}

To compute $\pi\left(10^{24}\right)$ with the zeros at our disposal we needed a sieve of width $\approx
6\times 10^{15}$ centred at $10^{24}$. We will discuss only locating the
primes in that interval, the prime powers being a trivial task by
comparison. 

Two basic methods were considered, sieving (necessarily segmented)
and a hybrid technique described by Galway \cite{Galway2004}. The latter proceeds by first eliminating all y-smooth numbers and then applying a base $2$ Fermat primality test to the remainder. Given a list of the (few) numbers in our range which are composite, y-rough and yet still pass the Fermat test, we are done. Our tests suggest that while an implementation of the Hybrid sieve would not be competitive at height $10^{24}$, the crossover might not be far away.

Our implementation used Atkin and Bernstein's sieve based on binary quadratic forms \cite{Atkin2004} to enumerate the sieving primes ($\leq
x^{1/2}$) which are then used in a segmented version of the sieve of Eratosthenes to delete composites in the target region.

For each sieve segment centred at $x_0$, we output
\begin{equation*}
\sum\limits_{p}1\textrm{ ,}\sum\limits_{p}(x_0-p)\textrm{ and }\sum\limits_{p}(x_0-p)^2.
\end{equation*}
By restricting the sizes of the segments, we can ensure that the entire
computation can be achieved using native $64$ bit integer instructions, with the exception of
the third sum which requires $128$ bit addition. However, this represents only
a small performance penalty on modern CPUs. These three
terms are then used to form an approximation to $\sum\limits_{p}\phi(p)$ by
Taylor series. In fact, three terms are not enough to give us the required
precision so we exploit the following lemma to derive a linear approximation
to the fourth (cubic) term as well.

\begin{lemma}\label{lem:cub_app}
If we approximate the real cubic $y=a_3x^3$ on the interval $x\in[-w,w]$ where $w>0$ with the line $y=ax$ with $a=\frac{3a_3w^2}{4}$, then the magnitude of the error over the interval is $\leq\frac{\left|a_3\right|w^3}{4}$. What is more, in terms of minimising the worst case error, this line is the best choice of any quadratic.
\begin{proof}
\begin{figure}[tbp]
\centering
\fbox{\includegraphics[width=0.7\linewidth] {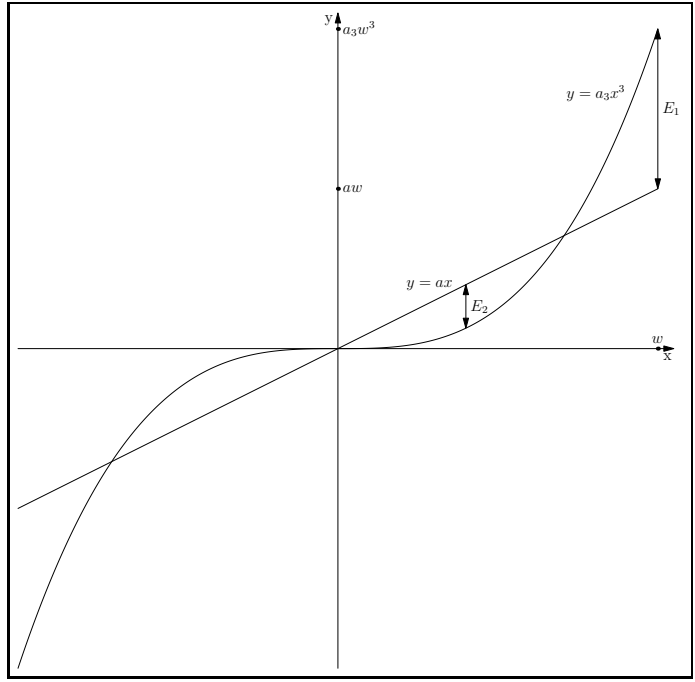}}
\caption{Approximating a cubic with a line (Lemma \ref{lem:cub_app})}
\label{fig:cubic}
\end{figure}

We refer to Figure \ref{fig:cubic}. Without loss of generality, take $a_3>0$. Since both $a_3x^3$ and $ax$ are odd, we consider only the interval $x\in[0,w]$. The error $E_1$ is simply $a_3w^3-aw$ and $E_2$ is at its maximum where the slopes of the line and the cubic are equal. This happens at $x=\sqrt{\frac{a}{3a_3}}$ so $E_2=\sqrt{\frac{a^3}{3a_3}}-\sqrt{\frac{a^3}{27a_3}}$. The worst case error follows from setting $E_1=E_2$ and solving for $a$.

The maximum error occurs $4$ times at $x\in\{\pm w,\pm \sqrt{\frac{a}{3a_3}}\}$. This means that any curve which improves on the line must be below the line at $x\in\{-w,\sqrt{\frac{a}{3a_3}}\}$ and above it at $x\in\{w,-\sqrt{\frac{a}{3a_3}}\}$. Thus, such a curve would have to cross the line at least $3$ times which is not possible for a quadratic.
\end{proof}
\end{lemma}

\section{Computing $\Phihat$}

The following lemmas give us a means of computing $\Re \Phihat(\rho)$.

\begin{lemma}
For $\Re s_0\neq0$ and $h\in\R$
\begin{equation}
\nonumber
\hat\phi(s_0+ih)=\hat\phi(s_0)\exp\left(ih(s_0\lambda^2+\log(x))\right)\frac{\exp\left(\frac{-\lambda^2h^2}{2}\right)}{1+\frac{ih}{s_0}}.
\end{equation}
\begin{proof}
We start with
\begin{equation}
\nonumber
\begin{aligned}
\hat\phi(s_0+ih)=\frac{\exp\left(\frac{\lambda^2(s_0+ih)^2}{2}\right)x^{s_0+ih}}{s_0+ih}
\end{aligned}
\end{equation}
and rearrange to get
\begin{equation*}
\begin{aligned}
\frac{\exp\left(\frac{\lambda^2s_0^2}{2}\right)x^{s_0}}{s_0}\frac{\exp\left(ih(s_0\lambda^2+\log(x))\right)\exp\left(\frac{-\lambda^2h^2}{2}\right)}{1+\frac{ih}{s_0}}.&\qedhere
\end{aligned}
\end{equation*}
\end{proof}
\end{lemma}

\begin{lemma}\label{lem:tay1}
Let $N\in2\N$, $\lambda,h>0$ and $\lambda h<1$. Then
\begin{equation}
\nonumber
\exp\left(\frac{-\lambda^2h^2}{2}\right)=\sum\limits_{n=0}^\frac{N}{2}\frac{(-1)^n}{n!}\left(\frac{\lambda^2h^2}{2}\right)^n+E_A,
\end{equation}
where
\begin{equation}
\nonumber
\left|E_A\right|\leq \frac{1}{\left(\frac{N}{2}\right)!}\left(\frac{\lambda^2h^2}{2}\right)^\frac{N}{2}.
\end{equation}
\begin{proof}
This function is entire and the restriction on $\lambda h$ makes the terms alternating in sign and decreasing.
\end{proof}
\end{lemma}

\begin{lemma}\label{lem:tay2}
Let $N\in\N$, $R=\left|\frac{h}{s_0}\right|$ and $|h|<|s_0|$. Then
\begin{equation}
\nonumber
\left(1+\frac{ih}{s_0}\right)^{-1}=\sum\limits_{n=0}^N\left(\frac{-ih}{s_0}\right)^n+E_B
\end{equation}
with
\begin{equation}
\nonumber
\left|E_B\right|\leq\frac{R^N}{1-R}.
\end{equation}
\begin{proof}
This function is analytic on the open disk $|h|<|s_0|$ and the missing terms form a geometric series.
\end{proof}
\end{lemma}

We can now fix some $N\in2\N$ and multiply these two (degree $N$) polynomials
to yield a single (degree $2N$) polynomial in $h$ which we can integrate
against $\exp(ih(\lambda^2+\log(x)))$ analytically.

We now start at $\Phihat\left(\frac{1}{2}\right)$ and move up the $\frac{1}{2}$ line in
small steps. We take the contribution from the highest-used $\rho$ to be zero and bound the error this approximation introduces.

\section{The Computation}

To obtain the $\bigoh\left(x^{1/2+\epsilon}\right)$ time complexity of Lagarias and
Odlyzko's algorithm, we should choose the free parameter $\lambda$ to equate the run
times of the sum over zeros and the sieving elements of the computation. In
fact, we biased the run time towards computing zeros, both to confirm RH holds to a height sufficient for this computation and since this data may have
application to future research.\footnote{To this end, Jonathan Bober has
  made about the first $36$ billion zeros available at \cite{Bober2012}.} We isolated all the zeros of $\zeta$ to a height of $30,610,046,000$ ($103,800,788,359$ zeros).  The technique used to locate these zeros is described in \cite{Platt2012} but is in essence a $\zeta$-specific, windowed version of Booker's algorithm from \cite{Booker2006}.

We set $\lambda=6273445730170391\times2^{-84}$ (note that this is exactly
representable in IEEE 754 double precision floating point). We used the first $69,778,732,700$ zeros to compute the sum (those to height $20,950,046,000$) which in turn dictated that we sieve a region of
width about $6\times 10^{15}$. As a
consequence, the
truncation error from summing over the zeros and from the sieve were together
$<0.989$.

With this choice of $\lambda$, we have
$\left|\Phihat\left(\frac{1}{2}\right)\right|<3\times 10^{13}$ so we need our
zeros to be located to an absolute accuracy of at least $25$ decimal
places.\footnote{Our zeros are located to an absolute accuracy of $\pm 2^{-102}$ which more than suffices.} Thus, we are forced to use multiple precision arithmetic, despite the
performance penalty this implies (up to a factor of $100$ compared with
hardware floating point).

As discussed earlier, we use interval arithmetic to manage the accumulation
of rounding errors during the computation and to this end we have extended
Revol and Rouillier's MPFI package \cite{Revol2002} in the obvious way to
handle complex arithmetic. Adopting interval arithmetic incurs another
performance penalty (a factor of about $3$ or $4$ this time).

The sieving (entirely in integer arithmetic) was performed on $352,800$ segments each
of width $2^{34}$ as dictated by memory constraints, further sub-divided to control the error in approximating $\phi$ by Taylor series. The actual computation of this Taylor approximation again requires multiple precision interval arithmetic.

The sum over zeros and the prime sieve all parallelise trivially and we used
the University of Bristol Bluecrystal Phase II cluster to perform all the
computations, consuming approximately $63,000$ CPU hours. In a personal communication Tom{\'a}s Oliveira e Silva indicated that computing $\pi\left(10^{23}\right)$ using the combinatorial method required about a month on a single computer. Assuming run time asymptotic to $\bigoh(x^{2/3})$ and $\bigoh(x^{1/2})$ for the combinatorial and analytic algorithms respectively, the crossover at which this implementation of the analytic algorithm would start to beat the combinatorial method would be in the region of $x=4\cdot 10^{31}$.
 
The result of the computation, after adding in the various error terms, was an interval straddling a single integer, so we have

\begin{theorem}
\begin{equation*}
\pi\left(10^{24}\right)=18,435,599,767,349,200,867,866.
\end{equation*}
\end{theorem}

We note that this agrees with the conditional result of B\"uthe, Franke, Jost and
Kleinjung.



\appendix
\section{Truncating the Sum over Zeros}\label{app:zero_trunc}

We require a rigorous bound for the error introduced by truncating the sum
over zeros in Theorem \ref{th:Gsum}. We proceed as follows.

For $t>0$, $t$ not the imaginary part of a zero of $\zeta$, define $N(t)$ to be the number of zeros of $\zeta(s)$ in the critical strip with $\Im s \in [0,t]$. 
\begin{lemma}\label{lem:Rosser}
Let $t\geq2$. Then
\begin{equation}
\nonumber
\left|N(t)-\frac{t}{2\pi}\log\left(\frac{t}{2\pi\e}\right)-\frac{7}{8}\right|<0.137\log(t)+0.443\log\log(t)+1.588.
\end{equation}
\begin{proof}
See \cite{Rosser1941} Theorem 19.
\end{proof}
\end{lemma}

\begin{lemma}\label{lem:Gbound}
For $x>1$, $T,\lambda>0$ and $\sigma\in[0,1]$ define
\begin{equation}
\nonumber
B(\sigma,T):=\exp\left(\frac{\lambda^2(1-T^2)}{2}\right)\left[\frac{x^\sigma}{T\log x}+\frac{1}{\lambda^2T^2x}\right].
\end{equation}
Then
\begin{equation}
\nonumber
\left|\Re \widehat{\Phi}\left(\sigma+iT\right)\right|\leq B(\sigma,T).
\end{equation}
\begin{proof}
We integrate along the contour running vertically down from $-1+i\infty$ to $-1+iT$, then right to $\sigma+iT$. For the horizontal contour we have

\begin{equation}
\nonumber
\begin{aligned}
\left|\int\limits^\sigma_{-1}\frac{\exp\left(\frac{\lambda^2(u+iT)^2}{2}\right)}{u+iT}x^{u+iT}\D u\right|&\leq\frac{\exp\left(\frac{\lambda^2(1-T^2)}{2}\right)}{T}\int\limits^\sigma_{-1}x^u\D u\\
&<\frac{\exp\left(\frac{\lambda^2(1-T^2)}{2}\right)}{T\log x}x^\sigma.
\end{aligned}
\end{equation}

For the vertical contour we have
\begin{equation}
\nonumber
\begin{aligned}
\left|\int\limits^T_\infty\frac{\exp\left(\frac{\lambda^2(-1+it)^2}{2}\right)}{-1+it}x^{-1+it}\D t\right|&\leq x^{-1}\exp\left(\frac{\lambda^2}{2}\right)\int\limits_T^\infty\frac{\exp\left(\frac{-\lambda^2t^2}{2}\right)}{t}\D t\\
&<\frac{\exp\left(\frac{\lambda^2}{2}\right)}{xT^2}\int\limits_T^\infty t\exp\left(\frac{-\lambda^2 t^2}{2}\right)\D t\\
&=\frac{\exp\left(\frac{\lambda^2(1-T^2)}{2}\right)}{\lambda^2T^2x}.&\qedhere
\end{aligned}
\end{equation}
\end{proof}
\end{lemma}

\begin{lemma}\label{lem:zerosum}
Let $T>0$, $\sigma\in[0,1]$ and $\alpha_T$ be such that $t^{\alpha_T} \geq N(t)$ for all
$t\geq T$. Then
\begin{equation}
\nonumber
\begin{aligned}
\sum\limits_{\Im \rho \geq T}B(\sigma,\Im \rho)\leq\exp\left(\frac{\lambda^2(1-T^2)}{2}\right)\left[\frac{x^\sigma}{T\log x}+\frac{1}{\lambda^2T^2x}\right]\left[\frac{\lambda^2T^2+2}{\lambda^2T^{2-\alpha_T}}-N(T)\right].
\end{aligned}
\end{equation}
\begin{proof}
Referring to Lemma \ref{lem:Gbound} and by writing $k_\sigma:=\exp\left(\frac{\lambda^2}{2}\right)\left[\frac{x^\sigma}{T\log x}+\frac{1}{\lambda^2T^2x}\right]$ we can majorise the sum with the Stieltjes integral
\begin{equation}
\nonumber
\int\limits_T^\infty k_\sigma\exp\left(\frac{-\lambda^2t^2}{2}\right) \D N(t).
\end{equation}
We now integrate by parts and majorise $N(t)$ with $t^{\alpha_T}$ to obtain
\begin{equation}
\nonumber
\begin{aligned}
\sum\limits_{\Im \rho >T}B(\sigma,\Im \rho)&\leq -k_\sigma\exp\left(\frac{-\lambda^2T^2}{2}\right)N(T)-\frac{k_\sigma}{T^{2-\alpha_T}}\int\limits_T^\infty\lambda^2t^3\exp\left(\frac{-\lambda^2t^2}{2}\right)\D t\\
&=k_\sigma\left[\frac{\lambda^2T^2+2}{\lambda^2T^{2-\alpha_T}}\exp\left(\frac{-\lambda^2 T^2}{2}\right)-\exp\left(\frac{-\lambda^2T^2}{2}\right)N(T)\right].&\qedhere
\end{aligned}
\end{equation}
\end{proof}
\end{lemma}
We note that the $\alpha_T$ referred to above can be computed using Lemma
\ref{lem:Rosser}.

We now consider the error introduced by truncating our sum over the zeros of $\zeta$. Let
$T_1$ be the height below which we find and use all the zeros, and $T_2$ the
height to which we know the Riemann Hypothesis holds.

\begin{lemma}\label{lem:E_1}
Let $E_1$ be real part of the error introduced by ignoring zeros with
imaginary part of absolute value $\in[T_1,T_2]$ (whose real parts are all
known to be $\frac{1}{2}$). Then 
\begin{equation}
\nonumber
\begin{aligned}
\left|E_1\right|\leq2\exp\left(\frac{\lambda^2(1-T_1^2)}{2}\right)\left[\frac{\sqrt{x}}{T_1\log x}+\frac{1}{\lambda^2T_1^2x}\right]\left[\frac{\lambda^2T_1^2+2}{\lambda^2T_1^{2-\alpha_{T_1}}}-N(T_1)\right].
\end{aligned}
\end{equation}
\begin{proof}
We apply Lemma \ref{lem:zerosum} with $\sigma=\frac{1}{2}$ and introduce a factor of $2$ for the zeros with negative imaginary part.
\end{proof}
\end{lemma}
We note this bound includes all the zeros with imaginary part $>T_2$ but their contribution will be negligible.

\begin{lemma}\label{lem:E_2}
Let $E_2$ be the real part of the error introduced by omitting the zeros with imaginary part $\not\in[-T_2,T_2]$ from the main
sum, (which do not necessarily have real part $=\frac{1}{2}$). Then
\begin{equation}
\nonumber
\left|E_2\right|\leq\exp\left(\frac{\lambda^2(1-T_2^2)}{2}\right)\left[\frac{x+1}{T_2\log x}+\frac{2}{\lambda^2T_2^2x}\right]\left[\frac{\lambda^2T_2^2+2}{\lambda^2T_2^{2-\alpha_{T_2}}}-N(T_2)\right].
\end{equation}
\begin{proof}
We pair each $\rho$ with $1-\overline{\rho}$ and take the worst case which is when one of the zeros has real part very close to $1$. The result then follows from Lemma \ref{lem:zerosum}.
\end{proof}
\end{lemma}

\iftoggle{longver}{
\section{Computing $f(t)$}\label{app:f}

In \cite{Odlyzko1988}, Odlyzko and Sch{\"o}nhage describe an efficient method to compute $\zeta$ at many equally spaced points simultaneously by recourse to the Fast Fourier Transform (FFT). In \cite{Booker2006}, Booker gives an algorithm for the rigorous computation of generic L-functions, again exploiting the efficiency of the FFT. Our desire for rigour draws us towards the latter of these, but its memory requirements would not be sustainable at large heights up the half line. We therefore developed a windowed version of Booker's algorithm for $\zeta$.

Recall we define $f$ as follows
\begin{equation}\label{eq:fdef}
f(t):=\pi^{-\frac{i(t+t_0)}{2}}\Gamma\left(\frac{\frac{1}{2}+i(t+t_0)}{2}\right)\exp\left(\frac{\pi(t+t_0)}{4}-\frac{t^2}{2h^2}\right)\zeta\left(\frac{1}{2}+i(t+t_0)\right).
\end{equation}

We proceed as follows:
\begin{enumerate}
\item Select $t_0,h>0$, $K\in\Z_{\geq 0}$ and $A,B>0$ such that $N=AB\in 2\N$. For
  $n=-\frac{N}{2}\ldots \frac{N}{2}-1$ and $k=0,1,\ldots,K$ compute
  $g\left(\frac{n}{A};k\right)$ where $g(t;k)$ is defined by
\begin{equation}
\nonumber
g(t;k):=\Gamma\left(\frac{\frac{1}{2}+i(t+t_0)}{2}\right)\exp\left(\frac{\pi
    (t+t_0)}{4}-\frac{t^2}{2h^2}\right)(-2\pi it)^k.
\end{equation}
\item By adding an appropriate error term, approximate
\begin{equation}
\nonumber
\widetilde{g}(n;k):=\sum\limits_{l\in\Z}g\left(\frac{n}{A}+lB;k\right).
\end{equation}
\item Use discrete Fourier transforms to compute
\begin{equation}
\nonumber
\widetilde{G}^{(k)}(m):=\sum\limits_{l\in\Z}G^{(k)}\left(\frac{m}{B}+lA\right),
\end{equation}
where
\begin{equation}
\nonumber
G(u):=\int\limits_{-\infty}^\infty g(t;0)\e(-t u) \D t.
\end{equation}
As customary, we define $\e(x):=\exp(2\pi i x)$.
\item Add an appropriate error term to recover $G^{(k)}\left(\frac{m}{B}\right)$ from $\widetilde{G}^{(k)}(m)$.
\item Use a series of convolutions to sum terms involving $G\left(\frac{m}{B}\right)$ and its derivatives yielding an approximation to $F\left(\frac{m}{B}\right)$, where
\begin{equation}
\nonumber
F(x):=\int\limits_{-\infty}^\infty f(t)\e(-tx)\D t.
\end{equation}
\item By adding an appropriate error term, approximate
\begin{equation}
\nonumber
\widetilde{F}(m):=\sum\limits_{l\in\Z}F\left(\frac{m}{B}+lA\right).
\end{equation}
\item Now use another discrete Fourier transform to compute
\begin{equation}
\nonumber
\widetilde{f}(n):=\sum\limits_{l\in\Z}f\left(\frac{n}{A}+lB\right).
\end{equation}
\item Finally, add another error term to recover $f\left(\frac{n}{A}\right)$ from $\widetilde{f}(n)$.
\end{enumerate}

Before presenting the necessary lemmas to make this outline concrete, we present the Theorem that allows us to bring the machinery of the FFT to bear.

\begin{theorem}\label{th:dft_pair}
Let $f$ be a function in the Schwartz space with Fourier transform $F$. Let
$N=AB$ with $A,B>0$ and define
$\tilde{f}(n):=\sum\limits_{l\in\Z}f\left(\frac{n}{A}+lB\right)$ and
$\tilde{F}(m):=\sum\limits_{l\in\Z}F\left(\frac{m}{B}+lA\right)$. Then, up to
a constant factor, $\tilde{f}$ and $\tilde{F}$ form a DFT pair of length $N$.
\begin{proof}
By Poisson summation we have
\begin{equation}
\nonumber
\begin{aligned}
\sum\limits_{l\in\Z}f(t+lB)&=\frac{1}{B}\sum\limits_{l\in\Z}F\left(\frac{l}{B}\right)\e\left(\frac{lt}{B}\right)\\
\tilde{f}(n)&=\frac{1}{B}\sum\limits_{l\in\Z}F\left(\frac{l}{B}\right)\e\left(\frac{ln}{N}\right).
\end{aligned}
\end{equation}
We now write $l=l^{'}N+m$ to get
\begin{equation}
\nonumber
\begin{aligned}
\tilde{f}(n)&=\frac{1}{B}\sum\limits_{m=0}^{N-1}\sum\limits_{l^{'}\in\Z}F\left(\frac{l^{'}N+m}{B}\right)\e\left(\frac{(l^{'}N+m)n}{N}\right)\\
&=\frac{1}{B}\sum\limits_{m=0}^{N-1}\e\left(\frac{mn}{N}\right)\tilde{F}(m).
\end{aligned}
\end{equation}
This is by definition an iDFT.
\end{proof}
\end{theorem}

The utility of this theorem will be apparent when $f$ and $F$ both decay quickly enough to allow $\tilde{f}(n)$ and $\tilde{F}(m)$ to be approximated by $f\left(\frac{n}{A}\right)$ and $F\left(\frac{m}{B}\right)$ respectively.

\begin{lemma}\label{lem:incgamint}
Define the incomplete Gamma function for $\Re s>0$ by \cite{abramowitz1964}
\begin{equation}
\nonumber
\Gamma(s,x):=\int\limits_x^\infty t^{s-1}\e^{-t}\D t.
\end{equation}
Then, given $\kappa>-1$ and $x,h>0$, we have
\begin{equation}
\nonumber
\int\limits_x^\infty w^\kappa\exp\left(\frac{-w^2}{2h^2}\right)\D w=2^{\frac{\kappa-1}{2}}h^{\kappa+1}\Gamma\left(\frac{\kappa+1}{2},\frac{x^2}{2h^2}\right).
\end{equation}
\begin{proof}
Substitute $t=\frac{w^2}{2h^2}$.
\end{proof}
\end{lemma}

\begin{lemma}\label{lem:gbound}
For $k\in\Z_{\geq 0}$, $t\in\R$ and $t_0,h>0$, recall that we define $g$ by
\begin{equation}
\nonumber
g(t;k):=\Gamma\left(\frac{\frac{1}{2}+i(t+t_0)}{2}\right)\exp\left(\frac{\pi
    (t+t_0)}{4}-\frac{t^2}{2h^2}\right)(-2\pi it)^k.
\end{equation}
Then
\begin{equation}
\nonumber
\left|g(t;k)\right|\leq4|2\pi t|^k\exp\left(\frac{-t^2}{2h^2}\right).
\end{equation}
\begin{proof}
We use the trivial bound $\left|\Gamma\left(\frac{\frac{1}{2}+ix}{2}\right)\e^\frac{\pi x}{4}\right|<4$.
\end{proof}
\end{lemma}

\begin{lemma}\label{lem:g_twid_err}
Let $n\in\left[-\frac{N}{2},\frac{N}{2}-1\right]$ and $B>h\sqrt{k}$. Then
\begin{equation}
\nonumber
\begin{aligned}
&\left|\sum\limits_{l\in\Znz}g\left(\frac{n}{A}+lB;k\right)\right|\leq\\
&8(\pi B)^k\left[\exp\left(\frac{-B^2}{8h^2}\right)+2^\frac{3k-1}{2}\left(\frac{h}{B}\right)^{k+1}\Gamma\left(\frac{k+1}{2},\frac{B^2}{8h^2}\right)\right].
\end{aligned}
\end{equation}
\begin{proof}
We consider the right tail from $n=-\frac{N}{2}$. The first term missing is
$g\left(\frac{B}{2};k\right)$ and $\frac{B}{2}$ is sufficiently large that our
bound for $g(t;k)$ (Lemma \ref{lem:gbound}) is decreasing. Thus we can split
off the first term and majorise the balance with an integral. This process
results in
\begin{equation}
\nonumber
\begin{aligned}
2\left[\left|g\left(\frac{B}{2};k\right)\right|+
\int\limits_1^\infty4(\pi B(2w-1))^k\exp\left(\frac{-(2w-1)^2B^2}{8h^2}\right)\D w\right].
\end{aligned}
\end{equation}
We now apply Lemma \ref{lem:incgamint} to the integral and the result follows.
\end{proof}
\end{lemma}

Thus, appealing to Lemma \ref{lem:g_twid_err} and choosing the parameters $B$ and $h$, we can control the
error introduced by using $g\left(\frac{n}{A};k\right)$ in place of $\widetilde{g}\left(n;k\right)$.

We now wish to compute values of $\widetilde{G}^{(k)}(m)$ from
$\widetilde{g}(n;k)$. The following lemma provides an efficient mechanism to
achieve this.

\begin{lemma}\label{lem:g_trans}
Up to a constant factor, the functions $\widetilde{g}(n;k)$ and $\widetilde{G}^{(k)}(m)$ form a discrete
Fourier transform pair of length $N$.
\begin{proof}
It is a standard result that given $f$ with Fourier transform $F$, the Fourier
transform of $x^kf(x)$ is $\left(\frac{i}{2\pi}\right)^k\frac{\D^k F(u)}{\D
  u^k}$. Thus $G^{(k)}(u)$ is the Fourier transform of $g(t;k)$ and the result
follows from Theorem \ref{th:dft_pair}.
\end{proof}
\end{lemma}

The result of the DFT above is $N$ values of $\widetilde{G}^{(k)}(m)$. The following
lemmas bound the error in using these values in place of $G\left(\frac{m}{B}\right)$.

\begin{lemma}\label{lem:C}
For $\sigma\in 2\Z_{>0} +1$ define
\begin{equation}
\nonumber
\begin{aligned}
&C(\sigma,t_0,h,k):=\\
&(2\pi)^k\int\limits_{-\infty}^\infty\left|\Gamma\left(\frac{\sigma+i(t+t_0)}{2}\right)\exp\left(\frac{\pi (t+t_0)}{4}-\frac{t^2}{2h^2}\right)\left(\frac{1}{2}+\sigma-it\right)^k\right|\D t.
\end{aligned}
\end{equation}
Then for $t_0>\sigma+\frac{1}{2}$ writing
\begin{equation}
\nonumber
X:=h^{m+1}2^\frac{\sigma-3-4k}{4}(2\sigma+1)^k\left(\Gamma\left(\frac{\sigma+1}{4}\right)-\Gamma\left(\frac{\sigma+1}{4},\frac{(2\sigma+1)^2}{8h^2}\right)\right)
\end{equation}
and
\begin{equation}
\nonumber
Y:=h^\frac{\sigma+1+2k}{2}2^\frac{\sigma-3+2k}{4}\Gamma\left(\frac{\sigma+1+2k}{4},\frac{(2\sigma+1)^2}{8h^2}\right)
\end{equation}
we have
\begin{equation}
\nonumber
\begin{aligned}
C(\sigma,t_0,h,k)\leq2^\frac{6k+7-\sigma}{4}\pi^\frac{2k+1}{2}\e^\frac{1}{2\sigma}\sum\limits_{l=0}^\frac{\sigma-1}{2}{\frac{\sigma-1}{2} \choose l}t_0^\frac{\sigma-2l-1}{2}
\left[X+Y\right].
\end{aligned}
\end{equation}

\begin{proof}
\begin{equation}
\nonumber
\begin{aligned}
&(2\pi)^k\int\limits_{-\infty}^\infty\left|\Gamma\left(\frac{\sigma+i(t+t_0)}{2}\right)\exp\left(\frac{\pi (t+t_0)}{4}-\frac{t^2}{2h^2}\right)\left(\frac{1}{2}+\sigma-it\right)^k\right| \D t\\
&\leq2^{k+1}\pi^k\int\limits_{0}^\infty\left|\Gamma\left(\frac{\sigma+i(t+t_0)}{2}\right)\exp\left(\frac{\pi (t+t_0)}{4}-\frac{t^2}{2h^2}\right)\left(\frac{1}{2}+\sigma-it\right)^k\right|\D t
\end{aligned}
\end{equation}
We now split the integral so that the upper bound becomes
\begin{equation}
\nonumber
\begin{aligned}
2^\frac{6k+7-\sigma}{4}\pi^\frac{2k+1}{2}\e^\frac{1}{2\sigma}&\left[\int\limits_0^{\sigma+\frac{1}{2}}(t+t_0)^\frac{\sigma-1}{2}(\sigma+1/2)^k\exp\left(\frac{-t^2}{2h^2}\right)\D t\right.\\
&\left.+\int\limits_{\sigma+\frac{1}{2}}^\infty(t+t_0)^\frac{\sigma-1}{2}t^k\exp\left(\frac{-t^2}{2h^2}\right)\D t\right].
\end{aligned}
\end{equation}
Now the first integral evaluates to
\begin{equation}
\nonumber
\sum\limits_{l=0}^\frac{\sigma-1}{2}{\frac{\sigma-1}{2} \choose l}t_0^\frac{\sigma-2l-1}{2}h^{l+1}2^\frac{l-1}{2}(\sigma+1/2)^k\left(\Gamma\left(\frac{l+1}{2}\right)-\Gamma\left(\frac{l+1}{2},\frac{(2\sigma+1)^2}{8h^2}\right)\right)
\end{equation}
and the second to
\begin{equation}
\nonumber
\sum\limits_{l=0}^\frac{\sigma-1}{2}{\frac{\sigma-1}{2} \choose l}t_0^\frac{\sigma-2l-1}{2}h^{l+k+1}2^\frac{l+k-1}{2}\Gamma\left(\frac{l+k+1}{2},\frac{(2\sigma+1)^2}{8h^2}\right).
\end{equation}
\end{proof}
\end{lemma}

\begin{lemma}\label{lem:Gkbound}
Let $\sigma\in 2\N +1$. Then $G^{(k)}(u)$ is bounded in absolute terms by
\begin{equation}
\nonumber
\begin{aligned}
&C(\sigma,t_0,h,k)\exp\left(\frac{(2\sigma+1)^2}{8h^2}-(2\sigma-1)\pi |u|\right)+\\
&2^{k+2}\pi^{k+1}\exp\left(\frac{-t_0^2}{2h^2}\right)\sum\limits_{l=0}^\frac{\sigma-1}{2} \frac{((2l+1/2)^2+t_0^2)^\frac{k}{2}}{l!}\exp\left(\frac{(4l+1)^2}{8h^2}-(4l+1)\pi|u|\right).
\end{aligned}
\end{equation}
\begin{proof}
First we consider $u\geq0$. We write
\begin{equation}
\nonumber
\left|G^{(k)}(u)\right|=\left|\int\limits_{-\infty}^\infty\Gamma\left(\frac{\frac{1}{2}+i(t+t_0)}{2}\right)\exp\left(\frac{\pi(t+t_0)}{4}-\frac{t^2}{2h^2}\right)(-2\pi it)^k\e(-tu)\D t\right|.
\end{equation}
Substituting $s=\frac{1}{2}+i(t+t_0)$, we now move the line of integration right to $\Re(s)=\sigma$ giving
\begin{equation}\label{eq:Gboundplus}
\begin{aligned}
&\left|G^{(k)}(u)\right|\leq\\
&\exp\left(\frac{(2\sigma-1)^2}{8h^2}-\pi u(2\sigma-1)\right)\times\\
&(2\pi)^k\int\limits_{-\infty}^\infty\left|\Gamma\left(\frac{\sigma+i(t+t_0)}{2}\right)\exp\left(\frac{\pi (t+t_0)}{4}\right)\exp\left(\frac{-t^2}{2h^2}\right)\left(\frac{1}{2}-\sigma-it\right)^k\right|\D t.
\end{aligned}
\end{equation}
For $u<0$, we move the line of integration left to $\Re(s)=-\sigma$, picking up the poles of $\Gamma\left(\frac{s}{2}\right)$ at $s=0,-2,\ldots,1-\sigma$. These give a contribution bounded by
\begin{equation}
\nonumber
2^{k+2}\pi^{k+1}\exp\left(\frac{-t_0^2}{2h^2}\right)\sum\limits_{l=0}^\frac{\sigma-1}{2}\frac{((2l+1/2)^2+t_0^2)^\frac{k}{2}}{l!}\exp\left(\frac{(4l+1)^2}{8h^2}+(4l+1)\pi u\right).
\end{equation}
The integral which remains is now
\begin{equation}
\nonumber
\begin{aligned}
&(2\pi)^k\exp\left(\frac{(2\sigma+1)^2}{8h^2}+(2\sigma+1)\pi u\right)\times\\
&\int\limits_{-\infty}^\infty \left|\Gamma\left(\frac{-\sigma+i(t+t_0)}{2}\right)\exp\left(\frac{\pi (t+t_0)}{4}-\frac{t^2}{2h^2}\right)(\sigma+\frac{1}{2}-it)^k\right| \D t.
\end{aligned}
\end{equation}
Finally, for our range of $\sigma$ and for $t\in\R$, we have
$|\Gamma(-\sigma/2+it)|<|\Gamma(\sigma/2+it)|$ and the result follows.
\end{proof}
\end{lemma}

We are free to chose a value of $\sigma$ that minimises this bound for a particular choice of $u$. We note that for $t_0$ large compared to $h$, $C(\sigma,t_0,h,k)$ is of order $t_0^\frac{\sigma-1+2k}{2}$.

\begin{lemma}\label{lem:G_twid_err}
For $m\in[0,N/2]$ and $\sigma\in2\N+1$ 
\begin{equation}
\nonumber
\begin{aligned}
&\Bigg|\sum\limits_{l\in\Znz}G^{(k)}\left(\frac{m}{B}+lA\right)\Bigg|\leq2^{k+3}\pi^{k+1} \exp\left(\frac{-t_0^2}{2h^2}\right)S+\\
&\;\;\;\;2\left(1+\frac{1}{A\pi(2\sigma-1)}\right)C(\sigma,t_0,h,k)\exp\left(\frac{(2\sigma+1)^2}{8h^2}-\frac{A\pi(2\sigma-1)}{2}\right)\\
\end{aligned}
\end{equation}
where $S$ is the sum
\begin{equation}
\nonumber
\sum\limits_{l=0}^\frac{\sigma-1}{2}\left(1+\frac{1}{A\pi(4l+1)}\right)\frac{((2l+1/2)^2+t_0^2)^\frac{k}{2}}{l!}\exp\left(\frac{(4l+1)^2}{8h^2}-\frac{A\pi(4l+1)}{2}\right).
\end{equation}
\begin{proof}
The left tail from $m=\frac{N}{2}$ majorises every case. The first term missing is $G^{(k)}\left(\frac{-A}{2}\right)$ and the remainder of the left tail is less in absolute terms than
\begin{equation}
\nonumber
\begin{aligned}
\int\limits_1^\infty&\Bigg[C(\sigma,t_0,h,k)\exp\left(\frac{(2\sigma+1)^2}{8h^2}-\frac{A\pi(2n-1)(2\sigma-1)}{2}\right)+\\
\nonumber
&2^{k+2}\pi^{k+1}\exp\left(\frac{-t_0^2}{2h^2}\right)\times\\
\nonumber
&\qquad\sum\limits_{l=0}^\frac{\sigma-1}{2}\frac{((2l+1/2)^2+t_0^2)^\frac{k}{2}}{l!}\exp\left(\frac{(4l+1)^2}{8h^2}-\frac{A\pi (2n-1)(4l+1)}{2}\right)\Bigg]\D n.
\end{aligned}
\end{equation}
\end{proof}
\end{lemma}
Now armed with values of $G^{(k)}\left(\frac{m}{B}\right)$ for several $k$, we wish to compute $F\left(\frac{m}{B}\right)$. We use the following result.
\begin{lemma}\label{lem:Fsum}
Let $F$ be the Fourier transform of $f$ (\ref{eq:fdef}). Then
\begin{eqnarray}
\nonumber
\Bigg|F(x)&-&\sum\limits_{j=1}^\infty\frac{1}{\sqrt{j}}\left(j\sqrt{\pi}\right)^{-it_0}G\left(x+\frac{\log\left(j\sqrt{\pi}\right)}{2\pi}\right)\Bigg|\\
\nonumber
&&=2\pi^{\frac{5}{4}}\exp\left(\frac{1}{8h^2}-\frac{t_0^2}{2h^2}-\pi x\right).
\end{eqnarray}
\begin{proof}
We start with $F(x)=\int\limits_{-\infty}^\infty f(t)\e(-tx)\D t$ and substitute $s=\frac{1}{2}+i(t+t_0)$. We then shift
the line of integration right to $\Re(s)=\sigma>1$ picking up the simple pole of
$\zeta(s)$ with residue $1$ at $s=1$. Now write $\zeta(s)$ as a sum (over $j$), interchange the sum and integral and move the line of integration back to $\frac{1}{2}$.
\end{proof}
\end{lemma}

The following lemma allows us to truncate the sum over $j$ at $J$.
\begin{lemma}\label{lem:f_hat_sum_err}
Let $x\geq 0$. Then
\begin{eqnarray}
\nonumber
&&\Bigg|\sum\limits_{j>J}\frac{1}{\sqrt{j}}\left(j\sqrt{\pi}\right)^{-it_0}G\left(x+\frac{\log\left(j\sqrt{\pi}\right)}{2\pi}\right)\Bigg|\\
\nonumber
&&\leq C(\sigma,t_0,h,0)\exp\left(\frac{(2\sigma-1)^2}{8h^2}\right)\pi^\frac{1-2\sigma}{4}\frac{J^{1-\sigma}}{\sigma-1}.
\end{eqnarray}
\begin{proof}
Take $x=0$ and apply (\ref{eq:Gboundplus}) of Lemma \ref{lem:Gkbound}.
\end{proof}
\end{lemma}

This suggests that $J$ will need to be in the region of $(t_0)^\frac{1}{2}$.

We need to calculate $F(x)$ on a lattice of points $u_m$ each $\frac{1}{B}$ apart. Using Taylor's Theorem with $K$ terms (see Lemma \ref{lem:tay_err} below for the truncation error) we can write
\begin{equation}\label{eq:tay}
F(x)\approx\sum\limits_{k=0}^{K-1}\sum_m\frac{G^{(k)}(x+u_m)}{k!}S_m^{(k)}
\end{equation}
where we set $\xi=\frac{1}{2B}$ and then
\begin{equation}
\nonumber
S_m^{(k)}:=\sum\limits_{\frac{\log\left(j\sqrt{\pi}\right)}{2\pi}\in[u_m-\xi,u_m+\xi)}\frac{1}{\sqrt{j}}\left(j\sqrt{\pi}\right)^{-it_0}\left(\frac{\log\left(j\sqrt{\pi}\right)}{2\pi}-u_m\right)^k.
\end{equation}
Now for each $k$, (\ref{eq:tay}) is a discrete convolution so computing our approximation to
$F\left(\frac{m}{B}\right)$ is achieved by summing the output of $K$
such convolutions. The following lemma provides the
bound on the error from truncating the Taylor series to $K$ terms.

\begin{lemma}\label{lem:tay_err}
Let $w\in[-\xi,\xi]$. Then we have
\begin{equation}
\nonumber
\left|\sum\limits_{k=K}^\infty\frac{G^{(k)}(u)w^k}{k!}\right|\leq\frac{2^{\frac{K+5}{2}}\pi^{K+\frac{1}{2}}h^{K+1}\xi^K}{\Gamma\left(\frac{K+2}{2}\right)}.
\end{equation}
\begin{proof}
\begin{equation}
\nonumber
\begin{aligned}
\left|\sum\limits_{k=K}^\infty\frac{G^{(k)}(u)w^k}{k!}\right|&\leq\sup_{u'\in(u-\xi,u+\xi]}\left|\frac{G^{(K)}(u')\xi^K}{K!}\right|\\
&\leq\sup_{u'\in(u-\xi,u+\xi]}\left|\int\limits_{-\infty}^\infty\frac{g(t;k)\xi^K\e(-u't)}{K!}\D t\right|\\
&\leq8\int\limits_0^\infty\frac{(2\pi t\xi)^K}{K!}\exp\left(\frac{-t^2}{2h^2}\right)\D t\\
&=\frac{2^\frac{3K+5}{2}\pi^K\xi^Kh^{K+1}\Gamma\left(\frac{K+1}{2}\right)}{\Gamma(K+1)}
\end{aligned}
\end{equation}
and the result follows from the duplication formula for $\Gamma$.
\end{proof}
\end{lemma}

Since this error term occurs $J$ times in (\ref{eq:tay}), weighted by $\frac{1}{\sqrt{j}}$ each time, we multiply it by $2\sqrt{J}-1$. 
\begin{lemma}\label{lem:f_hat_bound}
Let $\sigma\in 2\Z+1$ and $1<\sigma<t_0$. Then we have
\begin{equation}
\begin{aligned}
\nonumber
\left|F(x)\right|\leq&\zeta(\sigma)\pi^\frac{1-2\sigma}{4}C(\sigma,t_0,h,0)\exp\left(\frac{(2\sigma-1)^2}{8h^2}-\pi|x|(2\sigma-1)\right)\\
&+2\pi^\frac{5}{4}\exp\left(\frac{1}{8h^2}-\pi|x|-\frac{t_0^2}{2h^2}\right).
\end{aligned}
\end{equation}
\begin{proof}
Since $f(t)$ is real, its Fourier Transform $F(x)$ has the property $F(-x)=\overline{F(x)}$ so we need only consider $x\geq 0$. We write $s=\frac{1}{2}+i(t+t_0)$ and shift the line of integration right to $\Re(s)=\sigma$ encountering the pole of $\zeta(s)$ at $s=1$. This yields a residue smaller in absolute terms than
\begin{equation}
\nonumber
2\pi^\frac{5}{4}\exp\left(\frac{1}{8h^2}-\pi x-\frac{t_0^2}{2h^2}\right).
\end{equation} 
The remaining integral is then bounded in exactly the same fashion as in Lemma \ref{lem:Gkbound} using $|\zeta(\sigma+it)|\leq|\zeta(\sigma)|$ for $\sigma>1$ and $t\in\R$.
\end{proof}
\end{lemma}

\begin{lemma}\label{lem:f_hat_twid_err}
For $n\in[0,\frac{N}{2}]$ we have
\begin{equation}
\nonumber
\begin{aligned}
&\left|\sum\limits_{l\in\Znz}F\left(\frac{n}{B}+lA\right)\right|\leq\\
&2\zeta(\sigma)\pi^\frac{1-2\sigma}{4}C(\sigma,t_0,h,0)\exp\left(\frac{(2\sigma-1)^2}{8h^2}-\frac{A\pi(2\sigma-1)}{2}\right)\left(1+\frac{1}{A\pi(2\sigma-1)}\right)\\
&+4\pi^\frac{5}{4}\exp\left(\frac{1-4t_0^2}{8h^2}-\frac{\pi A}{2}\right)\left(1+\frac{1}{A\pi}\right).
\end{aligned}
\end{equation}
\begin{proof}
The left tail from $n=N/2$ majorises all other cases. The first term missing is $F\left(\frac{-A}{2}\right)$ and the remaining terms are bounded by 
\begin{equation}
\begin{aligned}
\nonumber
\int\limits_1^\infty&\left[\zeta(\sigma)\pi^\frac{1-2\sigma}{4}C(\sigma,t_0,h,0)\exp\left(\frac{-\pi A(2\sigma-1)(2n-1)}{2}\right)\right.\\
&\left.+2\pi^\frac{5}{4}\exp\left(\frac{1-4t_0^2}{8h^2}-\frac{\pi(2n-1)A}{2}\right)\right]\D n.
\end{aligned}
\end{equation}
\end{proof}
\end{lemma}

We now need to move from $\widetilde{F}$ to $\widetilde{f}$. The following
lemma provides the means:

\begin{lemma}
Up to a constant factor, the functions $\widetilde{f}$ and $\widetilde{F}$ form a discrete
Fourier transform pair of length $N$.
\begin{proof}
We defined $f$ and $F$ to be a Fourier transform pair and Theorem
\ref{th:dft_pair} therefore applies.
\end{proof}
\end{lemma}

Hence we can compute $N$ values of $\widetilde{f}(n)$ from our $N$
approximations to $\widetilde{F}(m)$ efficiently via a single DFT.

The final step is to extract approximations to $f\left(\frac{n}{A}\right)$ from our
values of $\widetilde{f}(n)$. The following two lemmas bound the error
introduced if we simply equate them.

\begin{lemma}\label{lem:f_bound}
Given $t\geq0$ and $t_0>\exp(\e)$, set $\beta=\frac{1}{6}+\frac{\log\log
  t_0}{\log t_0}$. Then
\begin{equation}
\nonumber
|f(t)|\leq12(t+t_0)^\beta\exp\left(-\frac{t^2}{2h^2}\right).
\end{equation}
\begin{proof}
By \cite{Cheng2004} we have for $t+t_0>\e$
\begin{equation}\label{eq:zeta_bound}
\left|\zeta\left(\frac{1}{2}+i(t+t_0)\right)\right|\leq 3(t+t_0)^\frac{1}{6}\log(t+t_0)
\end{equation}
so with $(t+t_0)>\exp(\e)$ and $\beta$ defined as above
\begin{equation}
\nonumber
(t+t_0)^\frac{1}{6}\log(t+t_0)\leq(t+t_0)^\beta.
\end{equation}
The factor of $4$ comes from the trivial bound for the Gamma factor.
\end{proof}
\end{lemma}

\begin{lemma}\label{lem:f_twid_err}
For $t\geq 0$ and $t_0>\exp(\e)$ set
$\beta=\frac{1}{6}+\frac{\log\log t_0}{\log t_0}$. Then providing $\frac{\beta
  h^2}{t_0}\leq \frac{B}{2}\leq t_0$ and $n\in\left[0,N-1\right]$ we have
\begin{equation}
\nonumber
\left|\sum\limits_{l\in\Znz}f\left(\frac{n-N/2}{A}+lB\right)\right|\leq 24(X+\frac{2^\beta h}{B}(Y+Z)),
\end{equation}
where
\begin{equation}
\nonumber
X=\left(\frac{B}{2}+t_0\right)^\beta\exp\left(-\frac{B^2}{8h^2}\right),
\end{equation}
\begin{equation}
\nonumber
Y=(t_0)^\beta\sqrt\frac{\pi}{2}\left(\erf\left(\frac{t_0}{h\sqrt{2}}\right)-\erf\left(\frac{B}{2h\sqrt{2}}\right)\right)
\end{equation}
and
\begin{equation}
\nonumber
Z=2^\frac{\beta-1}{2}h^\beta\Gamma\left(\frac{\beta+1}{2},\frac{B^2}{8h^2}\right).
\end{equation}
\begin{proof}
The lower bound on $B$ ensures that the bound of Lemma \ref{lem:f_bound} is decreasing for $t\geq\frac{B}{2}$. The worst case is when $n=0$ and for any $n$, the right tail majorises the left. The first missing term to the right is $f\left(\frac{B}{2}\right)$ and the remaining terms are majorised by
\begin{equation}
\nonumber
\begin{aligned}
&\left|12\int\limits_0^\infty\left(\frac{(2w+1)B}{2}+t_0\right)^\beta\exp\left(-\frac{\left(\frac{(2w+1)B}{2}\right)^2}{2h^2}\right)\D w\right|\leq\\
&\frac{12}{B}\left|\int\limits_\frac{B}{2}^{t_0}(2t_0)^\beta\exp\left(\frac{-t^2}{2h^2}\right)\D t +\int\limits_{t_0}^\infty(2t)^\beta\exp\left(\frac{-t^2}{2h^2}\right)\D t \right|.
\end{aligned}
\end{equation}
The result follows from Lemma \ref{lem:incgamint}.
\end{proof}
\end{lemma}}{}

\bibliographystyle{amsplain}
\bibliography{platt-pixv2.bbl}{}
\end{document}